\documentclass{amsart}

\usepackage{fancyhdr}
\usepackage{lastpage}
\usepackage{stmaryrd,yhmath}

\newcommand{\bdis}{\begin{displaymath}}
\newcommand{\edis}{\end{displaymath}}
\newcommand{\be}{\begin{equation}}
\newcommand{\ee}{\end{equation}}

\newcommand{\mbb}{\mathbb}
\newcommand{\mcal}{\mathcal}
\newcommand{\pd}{\partial}
\newcommand{\vp}{\varphi}

\newcommand{\zf}{\zeta\left(\frac{1}{2}+it\right)}

\newtheorem{lemma}[]{Lemma}

\theoremstyle{definition}

\newtheorem{cor}[]{Corollary}

\theoremstyle{remark}
\newtheorem{remark}[]{Remark}

\newtheorem*{mydef1}{{\bf Theorem}}

\newtheorem*{mydef61}{{\bf Formula 1}}

\newtheorem*{mydef62}{{\bf Formula 2}}

\numberwithin{equation}{section}

\begin{document}

\title{On the roots of the equation $Z'(t)=0$}

\author{Jan Moser}

\address{Department of Mathematical Analysis and Numerical Mathematics, Comenius University, Mlynska Dolina M105, 842 48 Bratislava, SLOVAKIA}

\email{jan.mozer@fmph.uniba.sk}


\begin{abstract}
We have proved in this paper that the Lindel\" of hypothesis generates essential contraction of distances between consecutive odd-order zeros of the function $Z'(t)$. This paper is the translation of the paper \cite{11}
into the English except part 8 that we added in order to point out the I. M. Vinogradov' scepticism on possibilities of the method of trigonometric sums. \\
\end{abstract}

\maketitle

\section{The result}

Let (see \cite{5}, pp. 79, 329)
\be \label{1.1}
\begin{split}
& Z(t)=e^{i\vartheta(t)}\zf, \\
& \vartheta(t)=-\frac t2\ln\pi+\text{Im}\ln\Gamma\left(\frac 14+\frac 12it\right)= \\
& = \frac t2\ln\frac{t}{2\pi}-\frac t2-\frac{\pi}{8}+\mcal{O}\left(\frac 1t\right).
\end{split}
\ee
In this paper we obtain a new consequence of the Lindel\" of hypothesis. Namely, the function $Z(t)$ oscillates in the interval
\bdis
(T,T+T^\tau)
\edis
where $\tau$ is arbitrary small positive fixed number. This implies that the number of maxima and minima of the function $Z(t)$ is unbounded (as $T\to\infty$). \\

Here is the survey of our results. Let
\bdis
S(a,b)=\sum_{0<a\leq n<b\leq 2a}n^{it},\quad b\leq\sqrt{\frac{t}{2\pi}}
\edis
be the elementary trigonometric sum (see \cite{8}, p. 31). The following theorem holds true.

\begin{mydef1}
If
\be \label{1.2}
|S(a,b)|<A(\Delta)\sqrt{a}t^\Delta,\quad 0<\Delta<\frac 16
\ee
then there is a root of odd order of the equation
\be \label{1.3}
Z'(t)=0
\ee
in the interval
\be \label{1.4}
(T,T+T^\Delta\psi(T))
\ee
where $\psi(T)$ is an arbitrarily slowly increasing function unbounded from above as $T\to\infty$ (as for example $\psi(T)=\ln\ln\dots\ln T$).
\end{mydef1}

\begin{remark}
Any root of the odd order of the equation (1.3) is the point of local maximum or local minimum of the function $Z(t)$.
\end{remark}

For example, in the case (see \cite{6})
\be \label{1.5}
\Delta=\frac{35}{216}+\epsilon,\quad \epsilon>0
\ee
we have the following

\begin{cor}
The interval
\be \label{1.6}
\left( T,T+T^{35/216+\epsilon}\right)
\ee
contains a point of maximum or minimum of the function $Z(t)$.
\end{cor}

On the Lindel\" of hypothesis we have (see \cite{4})
\be \label{1.7}
|S(a,b)|<A(\epsilon)\sqrt{a}t^\epsilon,
\ee
i. e. we have the following

\begin{cor}
On the Lindel\" of hypothesis the interval
\be \label{1.8}
\left( T,T+T^\epsilon\psi(T)\right)
\ee
contains a point of maximum or minimum of the function $Z(t)$.
\end{cor}

\begin{remark}
Thus the Lindel\" of hypothesis gives an 100\% improvement of the exponent $\frac{35}{316}$ in (1.6). Simultaneously, we have an 100\% improvement of all incoming exponents $\Delta$ of the kind (1.2), (1.5).
\end{remark}

Next, let $N_0'(T)$ denote the number of the zeros of the equation (1.3) for $t\in (0,T]$. We have the following

\begin{cor}
On the Lindel\" of hypothesis we have the estimate
\be \label{1.9}
N_0'(T+T^\tau)-N_0'(T)>A(\tau,\epsilon)T^{\tau-\epsilon},\ 0<\epsilon<\tau
\ee
where $\tau$ is an arbitrary small fixed number.
\end{cor}

There are the following reasons to study the roots of the equation (1.3):
\begin{itemize}
\item[(a)] the truth of the Riemann hypothesis itself is connected with the question on distribution of the roots of the equation (1.3) (see, for example, \cite{7}, p. 34, Corollary 3),
\item[(b)] if the Riemann hypothesis is true then the roots of the equation (1.3) are connected with the question: are there multiple ($\geq 2$) zeros of the function $\zf$?
\end{itemize}

\section{The main formulae}

The proof of our Theorem lies on the following

\begin{mydef61}
\be \label{2.1}
Z'(t)=-2\sum_{n\leq \sqrt{\frac{t}{2\pi}}}\frac{1}{\sqrt{n}}(\vartheta'-\ln n)\sin(\vartheta-t\ln n)+\mcal{O}(t^{-1/4}\ln t).
\ee
\end{mydef61}

The proof of this formula is situated in the parts 3 and 4 of this text. Next, in the part 5, the formula (2.1) is transformed into the following

\begin{mydef62}
\be \label{2.2}
\begin{split}
& Z'(t)=-2\sum_{n<P_0}\frac{1}{\sqrt{n}}\ln\frac{P_0}{n}\sin(\vartheta-t\ln n)+\mcal{O}(T^{-1/4}\ln T), \\
& t\in [T,T+H],\ H\in (0,\sqrt[4]{T}],\ P_0=\sqrt{\frac{T}{2\pi}}.
\end{split}
\ee
\end{mydef62}

Let the sequence $\{ \tilde{t}_\nu\}$ be defined by the formula
\be \label{2.3}
\vartheta(\tilde{t}_\nu)=\pi\nu+\frac{\pi}{2},\ \nu\in\mbb{N}.
\ee
We obtain in the part 6 the following statement.

\begin{lemma}
From (1.2) the estimate
\be \label{2.4}
\sum_{T\leq\tilde{t}_\nu\leq T+H}Z'(\tilde{t}_\nu)=\mcal{O}(T^\Delta\ln^2T)
\ee
follows.
\end{lemma}

Next, in the part 7, we obtain

\begin{lemma}
From (1.2) the formula
\be \label{2.5}
\sum_{T\leq\tilde{t}_\nu\leq T+H}(-1)^\nu Z'(\tilde{t}_\nu)=-\frac{1}{2\pi}H\ln^2\frac{T}{2\pi}+\mcal{O}(T^\Delta\ln^2T)
\ee
follows.
\end{lemma}

\begin{remark}
If
\bdis
T^\Delta=o(H)
\edis
then (2.5) is the asymptotic formula, i. e. if, for example,
\be \label{2.6}
H=T^\Delta\psi(T).
\ee
\end{remark}

Finally, from (2.4) and (2.5) we obtain the following

\begin{lemma}
\be \label{2.7}
\begin{split}
& \sum_{T\leq\tilde{t}_{2\nu}\leq T+H}Z'(\tilde{t}_{2\nu})=-\frac{1}{4\pi}H\ln^2\frac{T}{2\pi}+\mcal{O}(T^\Delta\ln^2T), \\
& \sum_{T\leq\tilde{t}_{2\nu+1}\leq T+H}Z'(\tilde{t}_{2\nu+1})=\frac{1}{4\pi}H\ln^2\frac{T}{2\pi}+\mcal{O}(T^\Delta\ln^2T).
\end{split}
\ee
\end{lemma}

In the case (2.6) the assertion of the Theorem follows from (2.7).

\section{Proof of the Formula 1 (the first part)}

\subsection{}

We use the formula (see \cite{15}, p. 72, $x=y=\sqrt{\frac{t}{2\pi}}$, $\eta=\sqrt{2\pi t}$, $m=[x]$)
\be \label{3.1}
\begin{split}
& \zeta(s)=\sum_{n=1}^m\frac{1}{n^s}+\chi(s)\sum_{n=1}^m\frac{1}{n^{1-s}}+\\
& +\frac{e^{-i\pi s}\Gamma(1-s)}{2\pi i}\left(\int_{C_1}+\int_{C_2}+\int_{C_3}+\int_{C_4}\right)\frac{w^{s-1}e^{-mw}}{e^w-1}{\rm d}w,\
s=\sigma+it
\end{split}
\ee
where $C_1,\dots,C_4$ are the segments binding the following points in the $w$-plane
\be \label{3.2}
\begin{split}
& \infty+i\eta(1+c),\ c\eta+i\eta(1+c), \\
& c\eta+i\eta(1+c),\ -c\eta+i\eta(1-c), \\
& -c\eta+i\eta(1-c),\ -c\eta-i(2m+1)\pi, \\
& -c\eta-i(2m+1)\pi,\ \infty-i(2m+1)\pi,
\end{split}
\ee
correspondingly, and
\bdis
0<c\leq \frac{1}{2}.
\edis
Putting
\bdis
s=\frac{1}{2}+it
\edis
into (3.1) and multiplying the last by $e^{i\vartheta(t)}$, (comp. \cite{15}, p. 79) we obtain
\be \label{3.3}
\begin{split}
& Z(t)=e^{i\vartheta(t)}\zf= \\
& =2\sum_{n\leq \alpha(t)}\frac{1}{\sqrt{n}}\cos(\vartheta-t\ln n)-\frac{1}{2\pi}e^{\pi t+i\vartheta}\Gamma\left(\frac{1}{2}-it\right)W(t)
\end{split}
\ee
where
\be \label{3.4}
W(t)=\left(\int_{C_1}+\int_{C_2}+\int_{C_3}+\int_{C_4}\right)\psi(t,w){\rm d}w=\int_{C(t)}\psi(t,w){\rm d}w,
\ee
and
\be \label{3.5}
\psi(t,w)=\frac{w^{-1/2+it}e^{-mw}}{e^w-1},\quad \alpha(t)=\sqrt{\frac{t}{2\pi}}.
\ee

\subsection{}
Putting
\be \label{3.6}
\Phi(t)=\sum_{n\leq \alpha(t)}\frac{1}{\sqrt{n}}\cos(\vartheta-t\ln n)
\ee
we obtain (let, for example, $\delta>0$)
\be \label{3.7}
\begin{split}
& \Phi(t+\delta)-\Phi(t)= \\
& =\sum_{n\leq \alpha(t)}\frac{1}{\sqrt{n}}\left[\cos(\vartheta(t+\delta)-(t+\delta)\ln n)-\cos(\vartheta(t)-t\ln n)\right]+ \\
& +\sum_{\alpha(t)<n\leq\alpha(t+\delta)}\frac{1}{\sqrt{n}}\cos(\vartheta(t+\delta)-(t+\delta)\ln n)=\Sigma_1+\Sigma_2.
\end{split}
\ee
Since
\be \label{3.8}
\alpha(t+\delta)-\alpha(t)=\mcal{O}\left(\frac{\delta}{\sqrt{t}}\right),
\ee
then we can choose $\delta$ in such a way that the interval
\bdis
(\alpha(t),\alpha(t+\delta)]
\edis
does not contain a natural number. Consequently,
\be \label{3.9}
\Sigma_2=0.
\ee

\subsection{}
Next, we have (see (3.4))
\be \label{3.10}
\begin{split}
& W(t+\delta)-W(t)=\int_{C(t+\delta)}\psi(t+\delta,w){\rm d}w-\int_{C(t)}\psi(t,w){\rm d}w= \\
& =\int_{C(t+\delta)}\psi(t+\delta,w){\rm d}w-\int_{C(t)}\psi(t+\delta,w){\rm d}w+ \\
& +\int_{C(t)}[\psi(t+\delta,w)-\psi(t,w)]{\rm d}w= \\
& =\int_{C(t+\delta)\cup\{-C(t)\}}\psi(t+\delta,w){\rm d}w+\int_{C(t)}[\psi(t+\delta,w)-\psi(t,w)]{\rm d}w=\\
& =\int_{C(t)}[\psi(t+\delta,w)-\psi(t,w)]{\rm d}w
\end{split}
\ee
by the Cauchy theorem, since
\bdis
\psi(t+\delta,w)
\edis
is the analytic function with respect to the variable $w$ (and $\delta$ is arbitrarily small) in the region bounded by the contour $C(t+\delta)\cup\{-C(t)\}$. Consequently,
\be \label{3.11}
W'(t)=\int_{C(t)}\frac{\partial \psi(t,w)}{\partial t}{\rm d}w.
\ee

\subsection{}
Hence, from (3.3), by (3.2), (3.3), we obtain
\be \label{3.12}
\begin{split}
& Z'(t)=-2\sum_{n\leq\alpha(t)}\frac{1}{\sqrt{n}}(\vartheta'-\ln n)\sin(\vartheta-t\ln n)+ \\
& +\frac{1}{2\pi i}\left\{\vartheta'-\frac{\Gamma'}{\Gamma}\left(\frac{1}{2}-it\right)-\pi i\right\}e^{i\vartheta+\pi t}\Gamma\left(\frac{1}{2}-it\right)W(t)- \\
& -\frac{1}{2\pi}e^{i\vartheta+\pi t}\Gamma\left(\frac{1}{2}-it\right)W'(t)= \\
& =Z_1(t)+Z_2(t)+Z_3(t).
\end{split}
\ee
We can apply the analysis from \cite{15}, pp. 71-74 in the case of the function
\bdis
e^{i\vartheta}e^{\pi t-i\pi/2}\Gamma\left(\frac{1}{2}-it\right)W(t)
\edis
and we obtain the estimate
\bdis
\mcal{O}(t^{-1/4}).
\edis
Furthermore, by \cite{15}, pp. 25, 221, we have
\be\label{3.13}
\frac{\Gamma'}{\Gamma}\left(\frac{1}{2}-it\right)=\mcal{O}(\ln t),\ \vartheta'(t)=\mcal{O}(\ln t),
\ee
and then we obtain
\be \label{3.14}
Z_2(t)=\mcal{O}(t^{-1/4}\ln t).
\ee

\section{Proof of the Formula 1 (the second part)}

Since by (3.5) we have
\be \label{4.1}
\frac{\pd\psi(t,w)}{\pd t}=i\ln w\frac{w^{-1/2+it}e^{-mw}}{e^w-1},
\ee
then we obtain the estimate of $Z_3(t)$ (see (3.12)) by the method \cite{15}, pp. 71-74 (we must to consider the influence of the factor $\ln w$).

\subsection{}
We have on the contour $C_4$ (see \cite{15}, p. 72)
\be \label{4.2}
\begin{split}
& w=u+i(2m+1)\pi, \\
& \ln w=\left\{\begin{array}{rcl} \mcal{O}(\ln u) & , & u>m, \\ \mcal{O}(\ln t) & , & u\in (-c\eta,m], \end{array} \right.
\end{split}
\ee
and, consequently,
\be \label{4.3}
\int_{C_4}=\mcal{O}\left\{\eta^{-1/2}e^{-5\pi t/4}\int_{-c\eta}^\infty e^{-mu}|\ln w|{\rm d}u\right\}=
\mcal{O}\left\{ e^{t(c-5\pi/4)}\ln t\right\}.
\ee

\subsection{}
On $C_3$ we have
\bdis
|\ln w|<A\ln t,
\edis
then (comp. \cite{15}, p. 72)
\be \label{4.4}
\int_{C_3}=\mcal{O}\left\{ e^{-t(\pi/2+A)}\ln t\right\}.
\ee

\subsection{}
On $C_1$ we have
\be \label{4.5}
|\ln w|=\left\{\begin{array}{rcl} \mcal{O}(\ln t) & , & 0<u\leq \pi\eta, \\ \mcal{O}(\ln u) & , & u>\pi\eta, \end{array} \right.
\ee
then (comp. \cite{15}, p. 72)
\be \label{4.6}
\begin{split}
& \int_{C_1}=\mcal{O}\left\{ \eta^{-1/2}\ln t\int_0^{\pi\eta} e^{-(\pi/2+A)t}{\rm d}u\right\}+
\mcal{O}\left\{\eta^{-1/2}\int_{\eta\pi}^\infty e^{-xu}\ln u{\rm d}u\right\}= \\
& = \mcal{O}\left\{\eta^{1/2}e^{-(\pi/2+A)t}\ln t\right\},\quad \eta=2\pi x.
\end{split}
\ee

\subsection{}
Since on $C_2$
\bdis
w=i\eta+\lambda e^{i\pi/4},\quad |\lambda|<\sqrt{2c\eta},\quad \ln w=\mcal{O}(\ln t),
\edis
then we obtain for the corresponding parts of the integral the following estimate (comp. \cite{15}, pp. 73, 74)
\be \label{4.7}
\begin{split}
& \mcal{O}(\eta^\sigma t^{-1/2}e^{-\pi t/2})\to \mcal{O}(\eta^{1/2}t^{-1/2}e^{-\pi t/2}\ln t), \\
& \mcal{O}(\eta^{\sigma-1}e^{-\pi t/2})\to \mcal{O}(\eta^{-1/2}e^{-\pi t/2}\ln t).
\end{split}
\ee

\subsection{}
Since (comp. \cite{15}, p. 74)
\be \label{4.8}
e^{i\vartheta+\pi t}\Gamma\left(\frac 12-it\right)=\mcal{O}(e^{\pi t/2}),
\ee
then by (3.12), (4.3), (4.4), (4.6), (4.7) we obtain the estimate
\be \label{4.9}
Z_3(t)=\mcal{O}(t^{-1/4}\ln t).
\ee
Finally, from (3.12) by (3.14), (4.9) we obtain (2.1).

\section{Proof of the Formula 2}

If
\bdis
H\in (0,\sqrt[4]{T}],\quad t\in [T,T+H],
\edis
then
\be \label{5.1}
\alpha(T+H)-\alpha(T)=\mcal{O}\left(\frac{H}{\sqrt{T}}\right)=\mcal{O}(T^{-1/4}),
\ee
and, by (3.13),
\be \label{5.2}
\begin{split}
& \sum_{\alpha(T)\leq n\leq\alpha(T+H)}\frac{1}{\sqrt{n}}(\vartheta'-\ln n)\sin(\vartheta-t\ln n)=\mcal{O}(T^{-1/4}\ln T).
\end{split}
\ee
Consequently, we have from (2.1) by (5.2)
\be \label{5.3}
Z'(t)=-2\sum_{n<P_0}\frac{1}{\sqrt{n}}(\vartheta'-\ln n)\sin(\vartheta-t\ln n)+\mcal{O}(T^{-1/4}\ln T)
\ee
where
\bdis
P_0=\alpha(T)=\sqrt{\frac{T}{2\pi}}.
\edis
Next, by the formulae (see \cite{15}, pp. 221)
\be \label{5.4}
\vartheta'(t)=\frac{1}{2}\ln\frac{t}{2\pi}+\mcal{O}\left(\frac 1t\right),\ \vartheta''(t)\sim\frac{1}{2t},
\ee
we have
\be \label{5.5}
\vartheta'(t)=\vartheta'(T)+\mcal{O}\left(\frac HT\right)=\ln P_0+\mcal{O}(T^{-3/4}).
\ee
Since
\be \label{5.6}
\sum_{n<P_0}\frac{1}{\sqrt{n}}\mcal{O}(T^{-3/4})=\mcal{O}(T^{-1/2}),
\ee
then we obtain from (5.3) by (5.5), (5.6)
\bdis
Z'(t)=-2\sum_{n<P_0}\frac{1}{\sqrt{n}}\ln\frac{P_0}{n}\sin(\vartheta-t\ln n)+\mcal{O}(T^{-1/4}\ln T),
\edis
i. e. the formula (2.2).

\section{Proof of the Lemma 1}

\subsection{}

Since (see (2.3))
\be \label{6.1}
\sin\{\vartheta(\tilde{t}_\nu)\}=(-1)^\nu,\quad  \cos\{\vartheta(\tilde{t}_\nu)\}=0,
\ee
then we obtain from (2.2) the formula
\be \label{2.2}
\begin{split}
& Z'(\tilde{t}_\nu)=2(-1)^{\nu+1}\sum_{n<P_0}\frac{1}{\sqrt{n}}\ln\frac{P_0}{n}\cos\{\tilde{t}_\nu\ln n\}+\mcal{O}(T^{-1/4}\ln T),
\end{split}
\ee
and for the distance
$$\tilde{t}_{\nu+1}-\tilde{t}_\nu$$
we have the formula
\be \label{6.3}
\tilde{t}_{\nu+1}-\tilde{t}_\nu=\frac{2\pi}{\ln\frac{T}{2\pi}}+\mcal{O}\left(\frac{H}{T\ln^2T}\right)=
\frac{\pi}{\ln P_0}+\mcal{O}\left(\frac{T^{-3/4}}{\ln T}\right)
\ee
(similarly to \cite{14}, p. 102; \cite{8}, (42)).

\subsection{}

Next, we have the following remarks:
\begin{itemize}
\item[(a)] the formula (6.3) for the difference $\tilde{t}_{\nu+1}-\tilde{t}_\nu$ is asymptotically equal to the formula for the difference
$t_{\nu+1}-t_\nu$ (comp. \cite{8}, (42)),
\item[(b)] the formula for $Z'(\tilde{t}_\nu)$ (see (6.2)) differs from the formula (see \cite{15}, p. 221)
\be \label{6.4}
Z(t_\nu)=2(-1)^\nu\sum_{n<P_0}\frac{1}{\sqrt{n}}\cos(t_\nu\ln n)+\mcal{O}(T^{-1/4})
\ee
by the factors
\bdis
-\ln\frac{P_0}{n},\quad \ln T
\edis
in the sum and in the error term, correspondingly.
\end{itemize}

\subsection{}

It is clear by (6.2) that our method of reduction
\bdis
\sum_{T\leq t_\nu\leq T+H}Z(t_\nu)\to W(t)
\edis
(see \cite{8}, (59)-(61)) is applicable also in the present case if we use the substitutions
\bdis
\frac{1}{\sqrt{n}}\to -\frac{1}{\sqrt{n}}\ln\frac{P_0}{n},\quad
\mcal{O}\{(\dots)\}\to\mcal{O}\{(\dots)\ln T\}.
\edis

\subsection{}

By (6.3) we have (see \cite{8}, (59)-(61))
\be \label{5.5}
\sum_{T\leq \tilde{t}_\nu\leq T+H}Z'(\tilde{t}_\nu)=-2\tilde{W}(T,H)+\mcal{O}(\ln^2T)
\ee
where ($n<P_0$)
\be \label{6.6}
\begin{split}
& \tilde{W}(T,H)=\frac{1}{2}(-1)^{\tilde{\nu}}\sum\frac{\ln\frac{P_0}{n}}{\sqrt{n}}\cos(\vp)+ \\
& +\frac{1}{2}(-1)^{N+\tilde{\nu}}\sum\frac{\ln\frac{P_0}{n}}{\sqrt{n}}\cos(\omega N+\vp)+ \\
& +\frac{1}{2}(-1)^{\tilde{\nu}}\sum\frac{\ln\frac{P_0}{n}\tan\frac{\omega}{2}}{\sqrt{n}}\sin(\vp)+ \\
& +\frac{1}{2}(-1)^{N+\tilde{\nu}+1}\sum\frac{\ln\frac{P_0}{n}\tan\frac{\omega}{2}}{\sqrt{n}}\sin(\omega N+\vp)+\mcal{O}(\ln^2T)= \\
& = \tilde{W}_1+\tilde{W}_2+\tilde{W}_3+\tilde{W}_4+\mcal{O}(\ln^2T),
\end{split}
\ee
and (see \cite{8}, (43), (50))
\be \label{6.7}
\begin{split}
& \tilde{t}_{\tilde{\nu}}=\min_{\tilde{t}_\nu\in [T,T+H]}\{\tilde{t}_\nu\},\quad
\tilde{t}_{\tilde{\nu}+N}=\max_{\tilde{t}_\nu\in [T,T+H]}\{\tilde{t}_\nu\}, \\
& \omega=\pi\frac{\ln n}{\ln P_0},\quad \vp=\tilde{t}_{\tilde{\nu}}\ln n.
\end{split}
\ee

\subsection{}

Since (see \cite{8}, (70))
\be \label{6.8}
\sum_{1\leq n<M\leq P_0}\frac{\cos\vp}{\sqrt{n}}=\mcal{O}(T^\Delta\ln T),
\ee
and the sequence
\bdis
\left\{\ln\frac{P_0}{n}\right\}
\edis
is decreasing and bounded by $A\ln T$ then we obtain by making use of  the Abel's transformation of the term $\tilde{W}_1$ (see (6.6) and (6.8))
\be \label{6.9}
\tilde{W}_1=\mcal{O}(T^\Delta\ln^2T),
\ee
and similarly
\be \label{6.10}
\tilde{W}_2=\mcal{O}(T^\Delta\ln^2T).
\ee

\subsection{}

For the sequence
\be \label{6.11}
\left\{\ln\frac{P_0}{n}\tan\frac{\omega}{2}\right\}
\ee
we have (see (6.7)) that
\be \label{6.12}
\begin{split}
& \tan\frac{\omega}{2}=\tan\left(\frac{\pi}{2}-\frac{\pi}{2}\frac{\ln\frac{P_0}{n}}{\ln P_0}\right)=
\cot\left(\frac{\pi}{2}\frac{\ln\frac{P_0}{n}}{\ln P_0}\right)=\cot X(n),
\end{split}
\ee
and consequently
\be \label{6.13}
\ln\frac{P_0}{n}\tan\frac{\omega}{2}=\frac{2}{\pi}\ln P_0 X(n)\cot X(n).
\ee
Next, the sequence
\bdis
\{ X(n)\cot X(n)\};\quad 0<X(n)\leq \frac{\pi}{2},\quad 1\leq n<P_0
\edis
is bounded by the value $1$, and it is increasing, since
\bdis
\frac{{\rm d}}{{\rm d}n}X\cot X=\frac{\pi}{2n\ln P_0}\frac{1-\frac{\sin 2X}{2X}}{\sin^2X}>0,\quad n\in [1,P_0).
\edis
Then, using the Abel's transformation and the estimate (comp. (6.8))
\bdis
\sum_{1\leq n<M\leq P_0}\frac{\sin\vp}{\sqrt{n}}=\mcal{O}(T^\Delta\ln T),
\edis
we obtain (see (6.13))
\be \label{6.14}
\sum_{n<P_0}X(n)\cot X(n)\frac{\sin\vp}{\sqrt{n}}=\mcal{O}(T^\Delta\ln T),
\ee
i. e. (see (6.6), (6.14))
\be \label{6.15}
\tilde{W}_3=\mcal{O}(T^\Delta\ln^2T),
\ee
and simultaneously
\be \label{6.16}
\tilde{W}_4=\mcal{O}(T^\Delta\ln^2T).
\ee
Finally, from (6.5) by (6.6), (6.9), (6.10), (6.15), (6.16) we obtain (2.4).

\section{Proof of the Lemma 2}

Let us remind that (see \cite{9}, (26))
\bdis
\sum_{T\leq t_\nu\leq T+H}(-1)^\nu Z(t_\nu)=\frac{1}{\pi}H\ln\frac{T}{2\pi}+\tilde{W}(T,H)+\mcal{O}\left(\frac{H^2}{T}\right),
\edis
where
\bdis
\tilde{W}=2\sum_{2\leq n<P_0}\frac{1}{\sqrt{n}}\sum_{T\leq t_\nu\leq T+H}\cos(t_\nu \ln n),
\edis
and (see \cite{9}, (51))
\be \label{7.1}
\tilde{W}=\mcal{O}(T^\Delta\ln T).
\ee
Since (see (6.3) and \cite{9}, (23))
\bdis
\sum_{T\leq \tilde{t}\nu\leq T+H}1=\frac{1}{2\pi}\ln\frac{T}{2\pi}+\mcal{O}(T^{-1/2})=\frac{1}{\pi}H\ln P_0+\mcal{O}(T^{-1/2}),
\edis
then we have (see (6.2))
\be \label{7.2}
\begin{split}
& \sum_{T\leq \tilde{t}\nu\leq T+H} (-1)^\nu Z'(\tilde{t}_\nu)= \\
& =-2\ln P_0\sum_{\tilde{t}_\nu}1-2\sum_{2\leq n<P_0}\frac{\ln\frac{P_0}{n}}{\sqrt{n}}\sum_{\tilde{t}_\nu}\cos(\tilde{t}_\nu\ln n)+\sum_{\tilde{t}_\nu}
\mcal{O}(T^{-1/4}\ln T)= \\
& =-2\ln P_0\left\{\frac{1}{\pi}H\ln P_0+\mcal{O}(T^{-1/2})\right\}-R(T,H)+\mcal{O}(\ln^2T),
\end{split}
\ee
and, of course, (comp. (7.1))
\bdis
\sum_{2\leq n<M\leq P_0}\frac{1}{\sqrt{n}}\sum_{T\leq \tilde{t}_\nu\leq T+H}\cos(\tilde{t}_\nu\ln n)=\mcal{O}(T^\Delta\ln T).
\edis
Hence, using the Abel's transformation on $R$, we obtain the estimate
\be \label{7.3}
R=\mcal{O}(T^\Delta\ln^2 T).
\ee
Finally, we obtain (2.5) from (7.2) by (7.3).

\appendix

\section{On I.M. Vinogradov' scepticism on possibilities of the method of trigonometric sums}

\subsection{}

I.M. Vinogradov analyzed in the Introduction to his monograph \cite{16} the possibilities of the method of trigonometric sums (Weyl's sums)
in the problem of estimation of the remainder term $R(N)$ in the asymptotic formula (see \cite{16}, p. 13)
\be \label{A.1}
\pi(N)-\int_2^N\frac{{\rm d}x}{\ln x}=R(N)
\ee
where $\pi(N)$ is the prime-counting function. Vinogradov made the following remark in this: \\

Obviously, it is very hard to make an essential progress in solution of the problem to find the order of the $R$-term (willing to find
$R=\mcal{O}(N^{1-c}),\ c=0.000001$) by making use of only some improvements of the H. Weyl's estimates and without making use of further important
progresses in the theory of the zeta-function.

\subsection{}

We will discuss in this section an analogue of the I.M. Vinogradov's scepticism in the case of estimation of the remainder term for the Hardy-Littlewood
integral in the formula
\be \label{A.2}
\int_0^T\left|\zeta\left(\frac{1}{2}+it\right)\right|{\rm d}t-T\ln T-(2c-1-\ln 2\pi)T=Q(T)
\ee
that is an analogue to (A.1). More generally, we ask whether there is a finer representation of the Hardy-Littlewood integral that the one
in the formula (A.2). \\

The first mathematician who applied the method of trigonometric sums to estimate $Q(T)$ was E.C. Titchmarsh (1934). He received the result (comp. \cite{15},
p. 123)
\bdis
Q(T)=\mcal{O}(T^{5/12+\epsilon}).
\edis
The crucial result in this field was obtained by Good (see \cite{1})
\be \label{A.3}
Q(T)=\Omega(T^{1/4}),\quad T\to\infty.
\ee
\begin{remark}
It follows from (A.3) that
\be\label{A.4}
Q(T)=\mcal{O}(T^{1/4+\epsilon}).
\ee
This Good's estimate still represents an unbounded and unremovable absolute error in the formula (A.2), and is the final eventuality for the method of
trigonometric sums in this question.
\end{remark}

Next, we have proved in our paper \cite{12} (90 years after the classical Hardy-Littlewood's paper \cite{2}) the following: there is an infinite set of
the almost exact representations of the Hardy-Littlewood integral
\be \label{A.5}
\int_0^T\left|\zeta\left(\frac{1}{2}+it\right)\right|{\rm d}t,
\ee
namely
\be \label{A.6}
\begin{split}
& \int_0^T\left|\zeta\left(\frac{1}{2}+it\right)\right|{\rm d}t=\vp_1(T)\ln\vp_1(T)+(c-\ln 2\pi)\vp_1(T)+c_0+\\
& +\mcal{O}\left(\frac{\ln T}{T}\right),\quad T\to\infty;\quad \vp_1(t)=\frac{1}{2}\vp(t),
\end{split}
\ee
(comp. \cite{13}, (9.1)), where $c$ is the Euler's constant, $c_0$ is the constant from the Titchmarsh-Kober-Atkinson formula (see \cite{15}, p. 141), and
$\vp(T)$ is a solution to the nonlinear integral equation
\bdis
\int_0^{\mu[x(T)]}Z^2(t)e^{-\frac{2}{x(T)}t}{\rm d}t=\int_0^T Z^2(t){\rm d}t;\quad \mu(y)\geq 7\ln y
\edis
in which each admissible function $\mu(y)$ (see \cite{12}) generates a solution
\bdis
y=\vp_\mu(T)=\vp(T).
\edis

\begin{remark}
The result (A.6) can be formulated as follows: the Jacob's ladders $\vp_1(T)$ (comp. (A.6) and the extension in \cite{12}, p. 415; $G[\vp(t)]$) are
the asymptotic solutions of the following new transcendental equation
\bdis
\int_0^T\left|\zeta\left(\frac{1}{2}+it\right)\right|{\rm d}t=V(T)\ln V(T)+(c-\ln 2\pi)V(T)+c_0.
\edis
\end{remark}

\begin{remark}
In the question on a possibility to find a finer representation of the Hardy-Littlewood integral (A.5) is the following situation:
\begin{itemize}
\item[(a)] we have for the classical representation (A.2)
\be \label{A.7}
\limsup_{T\to\infty}|Q(T)|=+\infty,
\ee
(see (A.3), (A.4)), i. e. we have the unbounded and unremovable absolute error term for all methods which use the estimation of trigonometric
sums,
\item[(b)] in our representation (A.6) of the Hardy-Littlewood integral that is absolutely independent on the method of trigonometric sums, we have
\be \label{A.8}
\lim_{T\to\infty}Q_1(T)=0;\quad Q_1(T)=\mcal{O}\left(\frac{\ln T}{T}\right),
\ee
i. e. we have a negligible error term,
\item[(c)] hence, the results (A.7), (A.8) confirm the validity of the analogue of I.M. Vinogradov's scepticism in the question of accuracy of the
representation of the Hardy-Littlewood integral (A.5).
\end{itemize}
\end{remark}

\subsection{}

We give in this part also other analogue of the I.M. Vinogradov's scepticism. First of all  we would like to describe our two main goals when working with
the Titchmarsh' sequence
$$\{ Z(t_\nu)\},$$
where $\{ t_\nu\}$ is the Gram's sequence. We wanted to:
\begin{itemize}
\item[(a)] improve the knowledge about the local variant of the classical Titchmarsh' formulae (see \cite{15}, pp. 221, 222)
\bdis
\begin{split}
& \sum_{\nu=\nu_0}^N Z(t_{2\nu})=2N+\mcal{O}(N^{3/4}\ln^{3/4}N), \\
& \sum_{\nu=\nu_0}^N Z(t_{2\nu+1})=-2N+\mcal{O}(N^{3/4}\ln^{3/4}N),
\end{split}
\edis
\item[(b)] prove a mean-value theorems for the function $Z(t)$ on related non-connected sets.
\end{itemize}

In order to solve the tasks (a) and (b) we have first used the method of trigonometric sums (see \cite{3}, pp. 260, 265; \cite{5}, p. 37). In the task
(a) we have improved the Titchmarsh' exponent as $\frac{3}{4}\to \frac{1}{6}$, i. e. we improved it by $77.7\%$. In the task (b) we obtained a new class of
mean-value theorems (see \cite{10}) corresponding to the exponent $\frac{1}{6}$
\bdis
\begin{split}
& \frac{1}{m\{ G_1(x)\}}\int_{G_1(x)}Z(t){\rm d}t\sim 2\frac{\sin x}{x}, \\
& \frac{1}{m\{ G_2(y)\}}\int_{G_2(y)}Z(t){\rm d}t\sim -2\frac{\sin y}{y}, \\
& 0<x,y\leq\frac{\pi}{2},\quad T\to\infty.
\end{split}
\edis
What concerns the task (a) -- we have essentially improved the Hardy-Littlewood exponent $\frac 14$ (since 1918) to the value $\frac 16$ (the
problem of estimation of the distance between neighboring zeros of the function $\zf$, see \cite{2}, p. 125, 177-184). Let us follow the sequence of
improvements of the mentioned Hardy-Littlewood exponent
\begin{itemize}
\item Moser -- 33\% improvement of the Hardy-Littlewood exponent $\frac 14$,
\item Karatsuba -- 6.25\% improvement of the exponent $\frac 16$,
\item Ivi\` c -- 0.19\% improvement of the Karatsuba exponent.
\end{itemize}

\begin{remark}
The sequence of improvements of the kind $6.25\%,\ 0.19\%,\dots$ shows that the scepticism of I.M. Vinogradov takes place also in the possibility
of successful application of the method of trigonometric sums in the problem of crucial improvement of the exponent $\frac{1}{6}$.
\end{remark}

\thanks{I would like to thank Michal Demetrian for helping me with the electronic version of this work.}


\begin{thebibliography}{29}
%
\bibitem{1}
A. Good, `Ein $\Omega$-Resultat f\" ur quadratische Mittel der Riemannschen Zetafunktion auf der kritische Linie`, Invent. Math. 41 (1977), 233-151.
%
\bibitem{2}
G.H. Hardy and J.E. Littlewood, `Contribution to the theory of the Riemann zeta-function and the theory of the distribution of primes`,
Acta Math. 41 (1918), 119-195.
%
\bibitem{3}
A. Ivic, `The Riemann zeta-function`, A Willey-Interscience Publications, New York, 1985.
%
\bibitem{4}
A.A. Karatsuba, `Basic analytic number theory`, Moscow, 1975 (in Russian).
%
\bibitem{5}
A.A. Karatsuba, `Complex analysis in number theory`, CRC Press, Boca Raton, Ann. Arbor, London, Tokyo, 1995.
%
\bibitem{6}
A. Kolesnik, `On the order of the Dirichlet $L$-function`, Pacific. J. Math. 82 (1979), 479-482.
%
\bibitem{7}
J. Moser, `Some properties of the Riemann zeta-function on the critical line`, Acta Arith. 26 (1974), 33-39, (in Russian).
%
\bibitem{8}
J. Moser, `On one sum in the theory of the Riemann zeta-function` Acta Arith., 31 (1976), 31-43; 40 (1981), 97-107, (in Russian).
%
\bibitem{9}
J. Moser, `One one theorem of Hardy-Littlewood in the theory of the Riemann zeta-function`, Acta Arith. 31 (1976), 45-51, Supplement in
Acta Arith. 35 (1979), 403-404, (in Russian).
%
\bibitem{10}
J. Moser, `New consequences of the Riemann-Siegel formula`, Acta Arith. 42 (1982), 1-10, (in Russian).
%
\bibitem{11}
J. Moser, `On the roots of the equation $Z'(t)=0$`, Acta Arith. 40 (1981), 97-107, (in Russian).
%
\bibitem{12}
J. Moser, `Jacob's ladders and the almost exact asymptotic representation of the Hardy-Littlewood integral`, Math. Notes, 88 (2010), 414-422,
arXiv: 0901.3973.
%
\bibitem{13}
J. Moser, `Jacob's ladders, the structure of the Hardy-Littlewood integral and some new class of nonlinear integral equations`, Proc. Stek. Inst.,
276 (2012), 208-221, arXiv: 1103.0359.
%
\bibitem{14}
E.C. Titchmarsh, `On van der Corput's method and the zeta-function of Riemann (IV)`, Quart. J. Math. 5 (1934), 98-105.
%
\bibitem{15}
E.C. Titchmarsh, `The theory of the Riemann zeta-function`, Clarendon Press, Oxford, 1951.
%
\bibitem{16}
I.M. Vinogradov, `The method of trigonometric sums in number theory`, NAUKA, Moscow, 2nd ed. 1980, (in Russian).
\end{thebibliography}
\end{document}